\documentclass[10pt]{article}
\usepackage[utf8]{inputenc}
\usepackage[english]{babel}
\usepackage{amsmath,amsthm,amsfonts}
\usepackage{graphicx}
\usepackage{bm}
\usepackage{tikz}
\usepackage{multirow}
\usepackage{color}

\graphicspath{ {./figs/} }

\newcommand{\grad}{\mathop{\rm grad}\nolimits}
\renewcommand{\div}{\mathop{\rm div}\nolimits}

\title{Constrained energy minimization based upscaling for coupled flow and mechanics}

\author{
Maria Vasilyeva \thanks{Institute for Scientific Computation, Texas A\&M University, College Station, TX 77843-3368 \& Department of Computational Technologies, North-Eastern Federal University, Yakutsk, Republic of Sakha (Yakutia), Russia, 677980. Email: {\tt vasilyevadotmdotv@gmail.com}.}
\and
Eric T. Chung \thanks{Department of Mathematics,
The Chinese University of Hong Kong (CUHK), Hong Kong SAR. Email: {\tt tschung@math.cuhk.edu.hk}.}
\and
Yalchin Efendiev \thanks{Department of Mathematics \& Institute for Scientific Computation (ISC),
Texas A\&M University,
College Station, Texas, USA. Email: {\tt efendiev@math.tamu.edu}.}
\and
Jihoon Kim
\thanks{Harold Vance Department of Petroleum Engineering, Texas A\&M University, College Station, Texas, USA. Email: {\tt jihoon.kim@tamu.edu}.}
}

\begin{document}

\maketitle

\begin{abstract}

In this paper, our aim is to present (1)
an embedded fracture model (EFM) for coupled flow and mechanics problem 
based on the dual continuum approach on the fine grid
and
(2) an upscaled model for the resulting fine grid equations.
The mathematical model is described by the coupled system of equation for displacement, 
fracture and matrix pressures.
For a fine grid approximation, we use the finite volume method for flow problem and 
finite element method for mechanics.
Due to the complexity of fractures, solutions have a variety of scales, 
and fine grid approximation results in a large discrete system. Our second focus in the 
construction of the upscaled coarse grid poroelasticity model for fractured media. 
Our upscaled approach is based on the nonlocal multicontinuum (NLMC) upscaling for 
coupled flow and mechanics problem, which involves computations of local basis functions
via an energy minimization principle. 
This concept allows a systematic upscaling for processes in the fractured porous media, 
and provides an effective coarse scale model whose degrees of freedoms have physical meaning.
We obtain a fast and accurate solver for the poroelasticity problem on a coarse grid and,
at the same time, derive a novel upscaled model.
We present numerical results for the two dimensional model problem.

\end{abstract}

\section*{Introduction}

In the reservoir simulation, mathematical modeling of the fluid flow and geomechanics in the 
fractured porous media plays an important role.
A coupled poroelastic models can help for better understanding of the processes in the 
fractured reservoirs.
In this work, we consider an embedded fracture model (EFM) for coupled flow and mechanics problems 
based on the dual continuum approach.
The mathematical model is described by the coupled system of equations for displacement
and fracture/matrix pressures \cite{salimzadeh2017three}.
Coupling of the fracture and matrix equations is derived from the mass exchange between 
the two continua (transfer term) and based on the embedded fracture model. 
For the geomechanical effect, we consider deformation of the porous matrix due to pressure change, 
where pressure plays a role of specific source term for deformation 
\cite{kim2010sequential, kim2011stability1, kim2011stability2, kolesov2014splitting, brown2016generalized, brown2016generalized2}.
Fundamentally, the system of equations is coupled between flow and geomechanics, where 
displacement equation includes the volume force, which is proportional to the pressure gradient, 
and the pressure equations include the term,
which describes the compressibility of the medium.

Fracture networks commonly have complex geometries with multiple scales, and usually have very 
small thickness compared to typical reservoir sizes.
Due to high permeability, fractures have a significant impact on the flow processes.
A common approach to the fracture modeling is to model them as lower dimensional problems 
\cite{martin2005modeling, d2012mixed, formaggia2014reduced, Quarteroni2008coupling}.
The result is a coupled mixed dimensional flow models, where we consider flow in the 
two domains (matrix and fracture) with mass transfer between them.
In this work, the fractures are not resolved by grid but included as an overlaying continuum 
with an exchange term between fracture and matrix that appears as an additional source 
(Embedded Fracture Model (EFM)) \cite{hkj12, ctene2016algebraic, tene2016multiscale}. 
This approach is related to the class of multicontinuum model 
\cite{barenblatt1960basic, warren1963behavior, douglas1990dual}. Instead of the dualcontinuum 
approach, we represent fractures directly using lower dimensional flow model embedded in a 
porous matrix domain.
In EFM, we have two independent grids for fracture networks and matrix, where simple structured meshes can be used for the matrix.

For geomechanics, we derive an embedded fracture model, where each fracture provides an additional 
source term for the displacement equation. This approach is based on the mechanics with dual 
porosity model \cite{wilson1982theory, zhao2006fully}.
In this model, we suppose displacement continuity on the fracture interface. For the 
discrete fracture model, a specific enrichment of the finite element space can be used for 
accurate solution of the elasticity problem with displacement discontinuity 
\cite{akkutlu2018multiscale}.
In this paper, we focus on the fully coupled poroelastic model for embedded fracture model and 
construct an upscaled model for fast coarse grid simulations.
For the fine grid approximation, we use the finite volume method (FVM) for 
flow problem and the finite element method (FEM) for geomechanics.
FVM is widely used as discretization for the simulation of flow problems 
\cite{bosma2017multiscale, ctene2017projection}. We use a cell centered finite volume 
approximation with two point flux approximation (TPFA) for pressure. FEM is typically used 
for approximating the solid deformation problem. We use a continuous Galerkin method with 
linear basis functions with accurate approximation of the coupling term.

Fine grid simulation of the processes in fractured porous media leads to very expensive simulations 
due to the extremely large degrees of freedoms. 
To reduce the cost of simulations, multiscale methods or upscaling techniques are used, for example, in \cite{houwu97, eh09, weinan2007heterogeneous, lunati2006multiscale, jenny2005adaptive}.
In our previous works, we presented multiscale model reduction techniques based on the Generalized multiscale finite element method (GMsFEM) for flow in fractured porous media \cite{akkutlu2015multiscale, chung2017coupling, efendiev2015hierarchical}. In GMsFEM approach, we solve a local spectral problem for the multiscale basis construction \cite{EGG_MultiscaleMOR, egh12, chung2016adaptive, CELV2015}. This gives us a systematic way to construct the missing degrees of freedom via multiscale
basis functions.
In this work, we construct an upscaled coarse grid poroelasticity model with 
embedded fracture model. 
Our approach uses the general concept of 
 nonlocal multicontinua (NLMC) upscaling for flow 
\cite{chung2017constraint, chung2017non}
and significantly generalized it to the coupled flow and mechanics problems.
The local problems for the upscaling
involves computations of local basis functions
via an energy minimization principle and the degrees of freedom are chosen such that
they represent physical parameters related to the coupled flow and mechanics problem.
We summarzie below the main goals of our work:
\begin{itemize}
\item a new fine grid embedded fracture model for poroelastic media (coupled system),
\item a new accurate and computationally effective fully coarse grid model for coupled 
multiphysics problem using NLMC whose degrees of freedoms have physical meaning on the coarse grid. 
\end{itemize}

Nonlocal multicontinua (NLMC) upscaling for processes in the fractured porous media provides 
an effective coarse scale model with physical meaning, and leads to a fast and accurate solver 
for coupled poroelasticity problem.
To capture fine scale processes at the coarse grid model, local multiscale basis functions 
are presented.
Constructing the basis functions based on the constrained energy minimization problem in the oversampled local domain is subject to the constraint that the local solution vanishes in other continua except the one for which it is formulated.
Multiscale basis functions have spatial decay property in local domains and separate background medium and fractures.
The proposed upscaled model has only one coarse degree of freedom (DOF) for each fracture network.
Numerical results show that our NLMC method for fractured porous media provides an accurate and efficient upscaled model on the coarse grid.

The paper is organized as follows.
In Section 1, we construct an embedded fracture model for poroelastic media.
Next, we construct fine grid approximation using FVM for flow problem and FEM for 
mechanics in Section 2.
In Section 3, we construct an upscaled coupled coarse grid poroelasticity model 
using NLMC method and present numerical results in Section 4.

\section{Embedded fracture model for poroelastic medium}

\begin{figure}[h!]
\centering
\includegraphics[width=0.7 \textwidth]{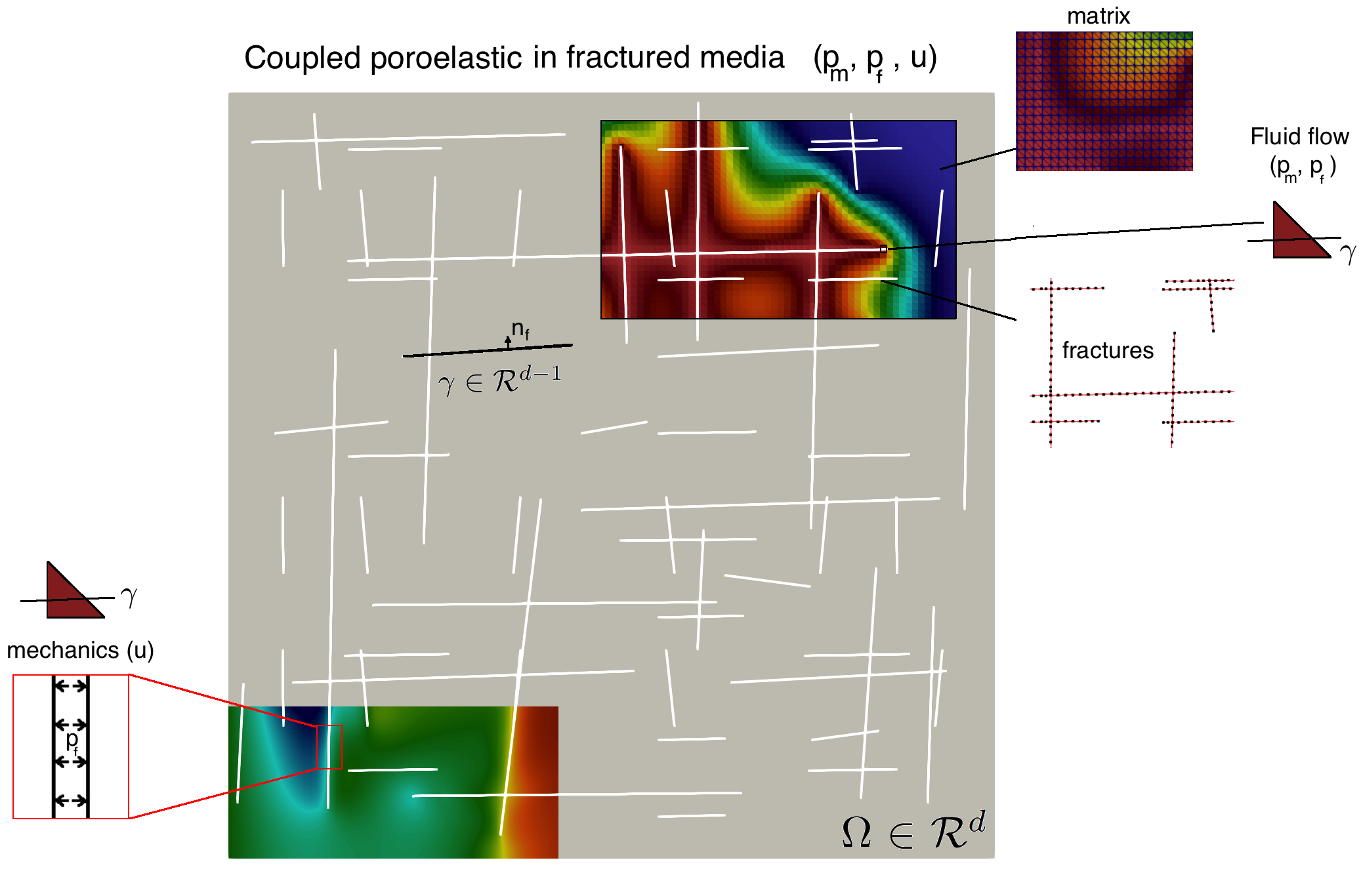}
\caption{Schematic illustration of the problem with embedded fracture model.}
\label{fig:sch}
\end{figure}

The proposed mathematical model of a coupled flow and mechanics in fractured poroelastic medium contains an interacting model for fluid flow in the porous matrix, flow in fracture network and mechanical deformation. 
The  matrix is assumed to be linear elastic and isotropic with now gravity effects. 
The mechanical and flow models are coupled through hydraulic loading on the fracture walls and using the effective stress concept \cite{wilson1982theory, salimzadeh2017three}. 
For fluid flow, we consider a mixed dimensional formulation, where we have a coupled problem for fluid flow in the porous matrix  in $\Omega \in \mathcal{R}^d$ (d = 2,3), and flow  in the fracture network on $\gamma \in \mathcal{R}^{d-1}$ (see Figure \ref{fig:sch} for $d =2$).

\textbf{Porous matrix flow model. } 
Using the mass conservation and Darcy law in the domain $\Omega$:
\begin{equation}
\label{eq:m-pe-2}
\frac{\partial m}{\partial t} + \div (\rho q_m) = \rho f_m, \quad 
q_m = - \frac{k_m}{ \nu_f} \grad p_m, \quad x \in \Omega,
\end{equation}
where $m$ is the fluid mass, $p_m$ is the matrix pressure, $q_m$ is Darcy velocity, $\nu_f$ is the viscosity, $\rho$ is the fluid density, and $f_m$ is the source term. 

Due to the motion of the solid skeleton and  Biot's theory, we have the following relationships
\cite{coussy2004poromechanics, kim2010sequential, kim2011stability1, kim2011stability2}
%
{
\begin{equation}
\label{eq:m2}
m - m_0  =
\rho \left( \frac{1}{M}  (p_m - p_0) + \alpha \varepsilon^v \right),
\end{equation}
where subscript $0$ means reference state, $\alpha$ is Biot coefficient, $M$ is Biot's modulus, $\varepsilon^v$ is the volumetric strain (the trace of the strain tensor, $\varepsilon^v = \text{tr } \varepsilon$) and  
\[
\frac{1}{M} = \Phi c_f + \frac{1}{N}, \quad 
\frac{1}{N} = \frac{\alpha - \Phi_0}{K_s}, \quad 
c_f = \frac{1}{\rho} \frac{d \rho}{d p_m}.
\]
Here $K_s$ is the solid grain stiffness, $c_f$ is the fluid compressibility and $\Phi$ is the Lagrange's porosity (also known as reservoir porosity). 

From equation \eqref{eq:m2}, we can express the reservoir porosity change induced by mechanical deformation as
\begin{equation}
\label{eq:m3}
\Phi - \Phi_0=\alpha \epsilon^v + \frac{1}{N} (p_m - p_0),
\end{equation}
}
The permeability of the matrix is updated using the current porosity by the power-law relationship 
\begin{equation}
\label{eq:m4}
k_m = k_0 \left(  \frac{\Phi}{\Phi_0} \right)^n.
\end{equation}
{
where the cubic law with $n = 3$ usually used \cite{wasaki2015permeability, zhang2008sorption}.
}

Therefore by assuming slightly compressible fluids, for the fluid flow in the porous matrix, we have the following parabolic equation
\[
\frac{1}{M} \frac{\partial p_m}{\partial t} 
+ \alpha \frac{\partial \varepsilon^v}{\partial t} 
- \div \left(  \frac{k_m}{ \nu_f} \grad p_m \right) 
=  f_m,
\]
defined in the domain $\Omega$.

{
For the case of fractures porous medium, we should add mass trasfer term between matrix and fracture
\begin{equation}
\label{eq:mm1}
\frac{1}{M} \frac{\partial p_m}{\partial t} 
+ \alpha \frac{\partial \varepsilon^v}{\partial t} 
- \div \left(  \frac{k_m}{ \nu_f} \grad p_m \right) + L_{mf}
=  f_m,
\end{equation}
where for the mass exchange between matrix and fracture, we assume a linear relationship
\[
L_{mf} = \beta_{mf} (p_m - p_f).
\]
This mass exchange term occurs only on the  fracture boundary.
}

\textbf{Fracture flow model. } 
For the highly permeable fractures, we use the following reduced dimension model for the fluid flow on $\gamma \in \mathcal{R}^{(d-1)}$ \cite{salimzadeh2017three, girault2016convergence}:
\begin{equation}
\label{eq:f1}
\frac{\partial (\rho b)}{\partial t} + \div (\rho \,  q_f) - \rho L_{mf} = \rho f_f, \quad 
q_f = - b \frac{k_f}{\nu_f} \grad p_f, \quad x \in \gamma,
\end{equation}
where $b$ is the fracture aperture, $p_f$ is the fracture pressure, $q_f$ is the average velocity of fluid along the fracture plane that can be calculated using the cubic low ($k_f = b^2$). For the calculation of the fracture aperture $b$, we can use following relation $b(t) = z p_f(t)$, where $z = \frac{2 (1-\nu^2)}{E}$ and deformation proportional to the fracture pressure $p_f$, {where $\eta$ is the Poisson's ratio, $E$ is the elastic modulus}  \cite{guo2012modeling, olorode17}. 

Since 
{
\begin{equation}
\label{eq:f2}
\frac{\partial (\rho b)}{\partial t} =
\rho \frac{\partial b}{\partial t} + b \frac{\partial \rho}{\partial t} =
\rho \left( \frac{\partial b}{\partial t}  + b c_f \frac{\partial p_f}{\partial t}\right),
\end{equation}
and by assuming slightly compressible fluids \cite{girault2016convergence}
\[
\rho \left(  \frac{\partial b}{\partial t}  + b c_f \frac{\partial p_f}{\partial t}\right) \approx
\rho_0 \left(  \frac{\partial b}{\partial t}  + b c_f \frac{\partial p_f}{\partial t}\right), 
\]\[
\div (\rho \, q_f) \approx \rho_0 \div q_f, \quad 
\rho L_{mf} \approx \rho_0 L_{mf}, \quad 
\rho f_f \approx \rho_0f_f.
\]
}
Therefore, we have the following equation on fracture $\gamma$
\begin{equation}
\label{eq:f3}
\begin{split}
\frac{\partial b}{\partial t}  + b c_f \frac{\partial p_f}{\partial t} - \div \left( b \frac{k_f}{\nu_f} \grad p_f \right) + L_{fm} = f_f, \quad x \in \gamma,
\end{split}
\end{equation}
where, for the mass exchange between matrix and fracture, we assume a linear relationship between the flux and pressure difference, namely, 
\[
L_{fm} = \beta_{fm} (p_f - p_m).
\] 

{Let 
$\beta_{fm} =  \eta_f \beta$ and $\beta_{mf} =  \eta_m \beta$, where $\beta$ is the trasnfer term proporsional to the matrix permeability, $\eta_f$ and $\eta_m$ are the geometric factors, that will be described in next section. }
Then we have 
\begin{equation}
\label{eq:mm2}
\begin{split}
& a_m \frac{ \partial p_m }{\partial t}  + \alpha \frac{\partial \varepsilon^v}{\partial t} 
- \div  (b_m \grad p_m) +  \eta_m \beta (p_m - p_f) =   f_m, 
\quad  x \in \Omega, \\
& a_f \frac{ \partial p_f }{\partial t} + \frac{\partial b}{\partial t}
- \div  (b_f \grad p_f) + \eta_f \beta (p_f - p_m) =   f_f.
\quad  x \in \gamma,
\end{split}
\end{equation}
where $a_m = 1/M$, $a_f = b c_f$, $b_m = k_m/{\nu_f}$, $b_f = b k_f/{\nu_f}$. 
{$\beta$ is the transfer term proportional to the matrix and fracture probabilities, $\eta_f$ and $\eta_m$ are geometric factors that will be defined in next section.} 
Pressure coupling term expresses the conservation of the flow rate (the fluid that is lost in the fractures goes into the porous matrix). 
Here we assume that the fractures have constant aperture, in general case, we can use $b(t) = z \, p_f(t)$ as a relationship between fracture width and pressure. 

\textbf{Mechanical deformation model. }
The balance of a linear momentum in the porous matrix is given by 
\begin{equation}
\label{eq:el1}
-\div \sigma_T = 0, \quad
\sigma_T = \sigma - \alpha p_m \mathcal{I}, \quad x \in \Omega,
\end{equation}
where 
$p_m$ is the matrix pressure, $\sigma_T$ is the total stress tensor, $\sigma$ is the effective stress \cite{salimzadeh2017three}.
Relation between the stress $\sigma$ and strain $\varepsilon$ tensors is given as
\[
\sigma = \lambda \varepsilon^v \mathcal{I} + 2 \mu \varepsilon (u), \quad
\varepsilon(u) = 0.5 (\nabla u + (\nabla u)^T),
\]
where {$u$ is the displacement vector in the porous matrix,} 
 and $\lambda$ and $\mu$ are the Lame's coefficients.

For incorporating of the fracture pressure into the model, we assume negligible shear traction on the fracture walls and consider normal tractions on the fractures \cite{salimzadeh2017three} with $\tau_f = -p_f n_f$, {where $n_f$ is the normal vector to the fracture surface}. 
After some manipulation, we obtain following equation in domain $\Omega$
\begin{equation}
\label{eq:mm3}
-\div \left( \sigma - \alpha p_m \mathcal{I} \right) {+} r_f p_f = 0,  \quad x \in \Omega,
\end{equation}
where $r_f$ comes from the integration over fracture surface ($\int_{\gamma} p_f \, n_f \,ds$) and contains direction of the fracture pressure influence. In presented model, we follow the classic dual porosity model and add fracture pressure effects as additional source (reaction) term. 
In more general case, fractures are modeled by an interface condition, where displacements have discontinuity across a fracture but stress is continuous \cite{akkutlu2018multiscale}.  

\section{Fine grid approximation of the coupled system}

Let $\mathcal{T}_h=  \cup_i \varsigma_i$ 
be a fine scale finite element partition of the domain $\Omega$ and $\mathcal{E}_{\gamma} = \cup_l \iota_l$ is the fracture mesh (see Figure \ref{fig:sch}).
The implementation is based on the open-source library FEniCS \cite{logg2009efficient, logg2012automated}. We use geometry objects for construction of the discrete system for coupled problem. For approximation of the flow part of the system, we use cell centered finite volume approximation with two point flux approximation. For displacement, we use Galerkin method with linear basis functions \cite{yoonspatial}.

In this work, we use the two dimensional problem for illustration of the robustness of our method.
In particular, we consider the following coupled system of equations for displacements (two displacements, $u_x$ and $u_y$) and fluid pressures (fracture and matrix, $p_f$ and $p_m$)
\begin{equation}
\label{eq:fd1}
\begin{split}
& a_m \frac{ \partial p_m }{\partial t}  + \alpha \frac{\partial \varepsilon^v}{\partial t} 
- \div  (b_m \grad p_m) +  \eta_m \beta (p_m - p_f) =   f_m, 
\quad  x \in \Omega, \\
& a_f \frac{ \partial p_f }{\partial t} + \frac{\partial b}{\partial t}
- \div  (b_f \grad p_f) + \eta_f \beta (p_f - p_m) =   f_f.
\quad  x \in \gamma, \\
& -\div \left( \sigma(u) - \alpha p_m \mathcal{I} \right) + r_f p_f = 0, \quad x \in \Omega.
\end{split}
\end{equation}

Using implicit scheme for approximation of time, a finite volume approximation for pressures and standard Galerkin method for displacements, we have following approximation
\begin{equation}
\label{eq:fd2}
\begin{split}
& \int_{\Omega} a_m \frac{ p_m - \check{p}_m }{\tau} d \Omega
+  \int_{\Omega} \alpha \frac{ \varepsilon^v - \check{\varepsilon}^v}{\tau}  d \Omega
 - \int_{\Omega}  \div  (b_m \grad p_m) d \Omega
 + \int_{\Omega}  \eta_m \beta (p_m - p_f )  d \Omega
 =  \int_{\Omega} f_m d \Omega, \\
& \int_{\gamma} a_f \frac{ p_f - \check{p}_f}{\tau} d \gamma
+ \int_{\gamma} \frac{b - \check{b}}{\tau} d \gamma
- \int_{\gamma} \div  (b_f \grad p_f) d \gamma
- \int_{\gamma} \eta_f \beta (p_m - p_f ) d \gamma
 =  \int_{\gamma} f_f  d \gamma, \\
& \int_{\Omega} (\sigma(u), \varepsilon(v)) d \Omega 
- \int_{\Omega} (\alpha p_m \mathcal{I}, \varepsilon(v)) d \Omega
+ \int_{\Omega} (r_f p_f, v)  d \Omega
= 0, \\
\end{split}
\end{equation}
where 
$(\check{p}_m, \check{p}_f, \check{u})$ are solutions from the previous times step and $\tau$ is the given time step.

Using the two point flux approximation for pressure equations, we obtain
\begin{equation}
\label{eq:fd3}
\begin{split}
& a_m \frac{ p_{m, i} - \check{p}_{m, i} }{\tau} |\varsigma_i| 
 +  \alpha \frac{ \varepsilon^v_i - \check{\varepsilon}^v_i}{\tau}  |\varsigma_i| 
 + \sum_j  T_{ij}  (p_{m, i} - p_{m, j})
 +  \beta_{il} (p_{m, i} - p_{f, l} )  
 =  f_m   |\varsigma_i|, \quad \forall i = 1, N^m_f \\
& a_f \frac{ p_{f, l} - \check{p}_{f, l}}{\tau}  |\iota_l| 
+ \frac{b_l- \check{b}_l}{\tau} |\iota_l| 
+ \sum_n W_{ln} (p_{f, l} - p_{f, n})
- \beta_{il} (p_{m, i} - p_{f, l} ) 
 =  f_f  |\iota_l|, \quad \forall l = 1, N^f_f
\end{split}
\end{equation}
where 
$T_{ij} = b_m |E_{ij}|/\Delta_{ij}$ ($|E_{ij}|$ is the length of {interface between cells} $\varsigma_i$ and $\varsigma_j$, $\Delta_{ij}$ is the distance between mid point of cells $\varsigma_i$ and $\varsigma_j$),  
$W_{ln} = b_f/\Delta_{ln}$ ($\Delta_{ln}$ is the distance between points $l$ and $n$), 
{ $|\varsigma_i|$ and $|\iota_l|$ is the volume of the cells  cells $\varsigma_i$ and $\iota_l$. }
{ $N^m_f$ is the number of cells in $\mathcal{T}_h$, 
$N^f_f$ is the number of cell for fracture mesh $\mathcal{E}_{\gamma}$.} 
Here, we use $\eta_m = 1 / |\varsigma_i|$ and $\eta_f = 1 / |\iota_l|$. Also,  $\beta_{il} = \beta$ if $\mathcal{E}_{\gamma} \cap \partial \varsigma_i = \iota_l$ and equals zero otherwise.

\textbf{Matrix form.} Combining the above schemes, we have following discrete system of equations for $y = (p_m, p_f, u_x, u_y)$ in the matrix form
\begin{equation}
\label{eq:fd4}
\left( \frac{1}{\tau} M + A \right) y = F, 
\end{equation}
where
\[
M = 
\begin{pmatrix}
M_m & 0 & 0  & 0 \\
0 & M_f & 0  & 0 \\
0 & 0 & 0  & 0\\
0 & 0 & 0  & 0
\end{pmatrix}, \quad
F =
\begin{pmatrix}
F_m 
+ \frac{1}{\tau} M_m \check{p}_m 
+ \frac{1}{\tau} (B_{m,x} + B_{m,y}) \check{u} \\
F_f + \frac{1}{\tau} M_f \check{p}_f \\
0\\
0
\end{pmatrix},
\]\[
{
A = 
\begin{pmatrix}
A_m + Q & -Q & \frac{1}{\tau} B_{m,x} & \frac{1}{\tau} B_{m,y} \\
-Q & A_f+Q & 0 & 0\\
-B_{m,x} & -B_{f,x} & D_{x} & D_{xy}\\
-B_{m,y} & -B_{f,y}	& D_{xy} & D_{y}
\end{pmatrix}, 
}
\]
\[
M_m = \{m^m_{ij}\}, \quad 
m^m_{ij} = 
\left\{\begin{matrix}
 a_m |\varsigma_i| / \tau & i = j, \\ 
0 & i \neq j
\end{matrix}\right. , \quad 
M_f = \{m^f_{ln}\}, \quad 
m^f_{ln} = 
\left\{\begin{matrix}
 a_f |\iota_l| / \tau & l = n, \\ 
0 & l \neq n
\end{matrix}\right. ,
\]\[
A_m = \{T_{ij}\}, \quad 
A_f = \{W_{ln}\}, \quad 
Q = \{q_{il}\}, \quad 
q_{il} = 
\left\{\begin{matrix}
\beta & i = l, \\ 
0 & i \neq l
\end{matrix}\right. ,
\]\[
F_m = \{f^m_i\}, \quad f^m_i = f_m |\varsigma_i|, \quad
F_f = \{f^f_l\}, \quad f^m_i = f_f |\iota_l|.
\]
Here $D$ is the elasticity stiffness matrix
{
\[
D_{x} = [d^{x}_{ij}] = 
\int_{\Omega}  \sigma_x(\psi_i) : \varepsilon_x(\psi_j) \, d\Omega, \quad 
D_{y} = [d^{y}_{ij}] = 
\int_{\Omega}  \sigma_y(\psi_i) : \varepsilon_y(\psi_j) \, d\Omega, \quad 
D_{xy} = [d^{xy}_{ij}] = 
\int_{\Omega}  \sigma_x(\psi_i) : \varepsilon_y(\psi_j) \, d\Omega, 
\]\[
B_{m,x} = [b^{m,x}_{ij}] = 
\int_{\Omega} (\alpha p_{m,i}, \varepsilon_x(\psi_j)) d \Omega, \quad 
B_{m,y} = [b^{m,y}_{ij}] = 
\int_{\Omega} (\alpha p_{m,i}, \varepsilon_y(\psi_j)) d \Omega, 
\]\[
B_{f, x} = [b^{f,x}_{lj}] = 
- \int_{\Omega} (r_f p_{f,l}, \psi_j)  d \Omega, \quad
B_{f, y} = [b^{f,y}_{lj}] = 
- \int_{\Omega} (r_f p_{f,l}, \psi_j)  d \Omega,
\]  
with linear basis functions $\psi_i$ and 
$\sigma = 
\begin{pmatrix} 
\sigma_{x} &  \sigma_{xy} \\
\sigma_{yx} &  \sigma_{y} \\
\end{pmatrix}$ and 
$\varepsilon = 
\begin{pmatrix} 
\varepsilon_{x} &  \varepsilon_{xy} \\
\varepsilon_{yx} &  \varepsilon_{y} \\
\end{pmatrix}$ 
}

We remark that the dimension of fine grid problem is given by
\[
N_f = N^m_f + N^f_f + 2 N^v_f,
\]
where $N^v_f$ is the number of vertices on the  fine grid.

\section{Coarse grid upscaled model for coupled problem}

Consider a coarse grid partition $\mathcal{T}_H = \{ K_i\}$ of the domain, where $K_i$ is the $i$-th coarse cell. 
Let $K^+_i$ be the oversampled region for the coarse cell $K_i$ obtained by enlarging $K_i$ by a few coarse grid blocks. 
For our coarse grid approximation, we will construct multiscale basis functions using the nonlocal multicontinua method (NLMC)\cite{chung2017non}. 
In general, the construction of the multiscale basis functions starts with an auxiliary  space,
which is constructed  by solving local spectral problems \cite{chung2017constraint}, and then we take eigenvectors that correspond to small (contrast dependent) eigenvalues as basis functions. These spectral basis functions
represent the channels (high contrast features). 
Using the auxiliary space, the target multiscale space is obtained by solutions of constraint energy minimization problems in oversampling domain $K^+_i$.
subject to { a set of orthogonality conditions} related to the auxiliary space. 
{ More precisely, for each auxiliary basis function, we will find a corresponding multiscale basis function such that 
it is orthogonal to all other auxiliary basis functions 
with respect to a weight inner product.
Our} basis functions have a nice decay property away from the target coarse element. 
In this paper, we use the NLMC method. 
In the NLMC  method, we use a simplified construction that separate continua in each local domain $K_i$ (coarse cell). 
Instead of using an auxiliary space, we obtain the required basis functions by minimizing an energy over an oversampling domain $K_i^+$
subject to the conditions that the minimizer has mean value zero on all fractures and matrix except the fracture or matrix
that the basis function is formulated for. 
The resulting
multiscale basis functions have a spatial decay property in local domains and separate background medium and fractures.

\begin{figure}[h!]
\centering
\includegraphics[width=1.0 \textwidth]{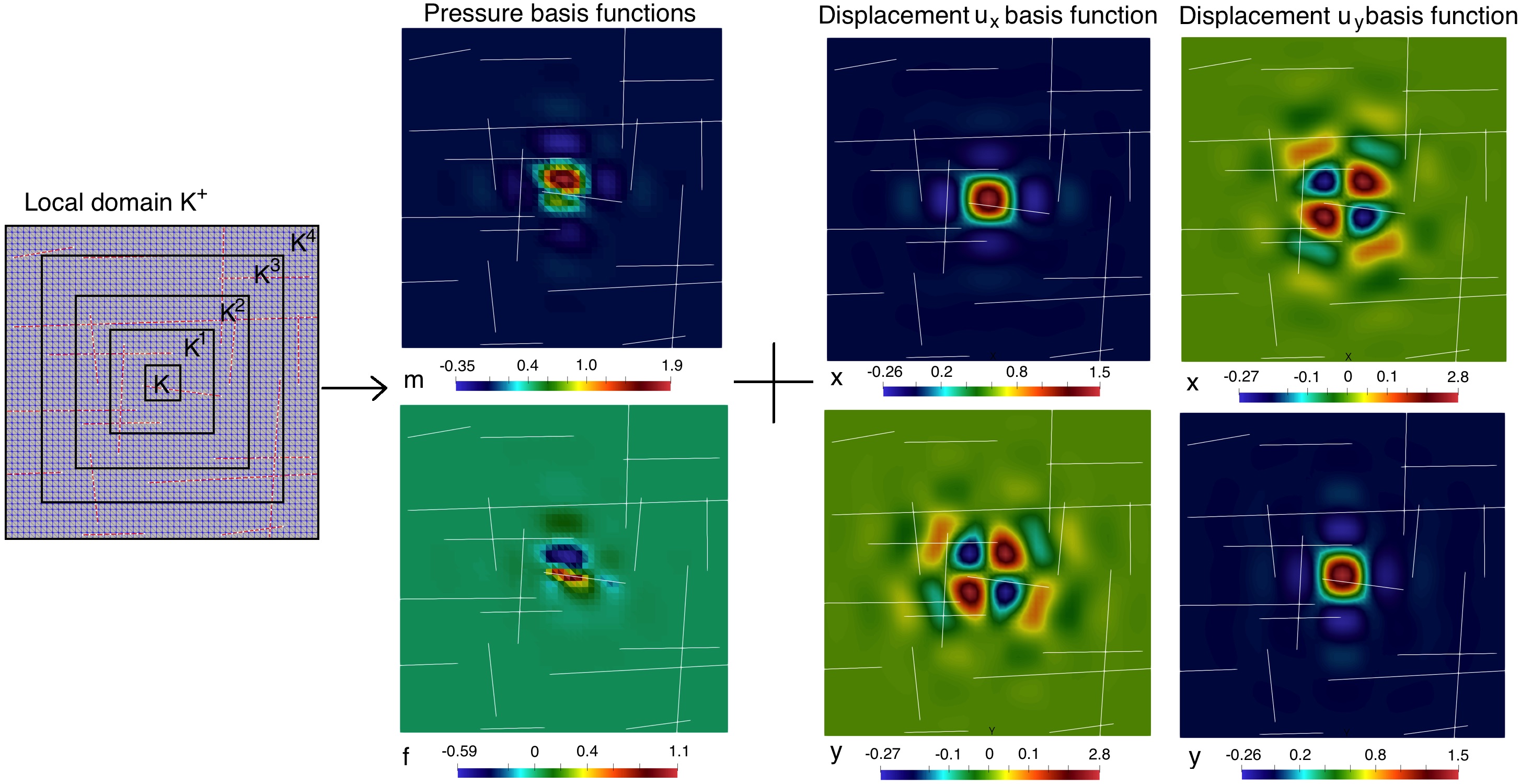}
\caption{Multiscale basis functions on mesh $20 \times 20$ for local domain $K^4$ for pressures and displacements.}
\label{fig:msbf20}
\end{figure}

For the fractures, we write $\gamma = \cup_{l = 1}^L \gamma^{(l)}$, where $\gamma^{(l)}$ is the $l$-th fracture network and $L$ is the total number of fracture networks. 
We also write $\gamma_j = \cup_{l = 1}^{L_j} \gamma_j^{(l)}$, where $\gamma^{(l)}_j = K_j \cap \gamma^{(l)}$ is the fracture  inside the coarse cell $K_j$ and  $L_j$ is the number of fractures in $K_j$. 
Again, for construction of multiscale basis functions, we solve constrained energy minimization problem in the oversampled local domains subject to the constraint that the local solution has zero mean on other continua except the one for which it is formulated. 

For  construction for coarse grid approximation for the coupled problem, we construct multiscale basis function for $(p_m, p_f, u_x, u_y)$. 
For simplicity, we ignore the coupling term between pressure and  displacements and find multiscale basis functions for pressure and displacements separately. In general coupled poroelastic basis functions can be constructed using the coupled poroelastic system and the constrained energy minimization principle. 

\textbf{Multiscale basis function for matrix and fracture pressures.} 
To define our multiscale basis functions, we will minimize an energy subject to some constraints. In the following, we will define the constraints. 
We remark that we will find a set of multiscale basis functions for each coarse cell $K_i$, and these basis functions have support in $K_i^+$. Thus
the following constraints are needed for each $K_i$, and they are defined within $K_i^+$.
For each coarse cell $K_j \in K_i^+$: \\
(1) background medium ($\psi^0_{l}$) :
\[
\int_{K_j} \psi_0^i \, dx = \delta_{i,j}, \quad 
\int_{\gamma^{(l)}_j} \psi_0^i \, ds = 0, \quad  l=\overline{1, L_j},
\]
We note that these constraints are defined for the matrix part in $K_i$, and they require the resulting basis function to have mean value on each continuum in $K_i^+$
except the continuum corresponding to the matrix part in $K_i$. \\
(2) $l$-th fracture network in $K_i$ ($\psi_l^i$):
\[
\int_{K_j} \psi_l^i \, dx = 0, \quad 
\int_{\gamma^{(l)}_j} \psi_l^i \, ds = \delta_{i,j}\delta_{m,l}, \quad  l=\overline{1, L_j},
\]
where $L_i$ is the number of fracture networks in $K_i$. 
We note that these constraints are defined for a fracture network in $K_i$, and they require the resulting basis function to have mean value on each continuum in $K_i^+$
except the continuum corresponding to a specific fracture network in $K_i$

For the construction of the multiscale basis functions, 
{we solve the following local problems in $K_i^+$ using an operator restricted in $K_i^+$. }
This results in solving the following local problems in $K_i^+$: 
\begin{equation}
\label{eq:basis}
{
\begin{pmatrix}
A_m^{K^+_i}+ Q^{K^+_i} & -Q^{K^+_i} & C^T_m & 0 \\
-Q^{K^+_i} & A_f^{K^+_i} + Q^{K^+_i} & 0 &  C^T_f \\
C_m & 0 & 0 & 0 \\
0 & C_f & 0 & 0 \\
\end{pmatrix} 
}
\begin{pmatrix}
\psi_m \\
\psi_f \\
\mu_m \\
\mu_f \\
\end{pmatrix} = 
\begin{pmatrix}
0 \\
0 \\
F_m \\
F_f \\
\end{pmatrix}
\end{equation}
with  zero Dirichlet boundary conditions on $\partial K^+_i$ for $\psi_m$ and $\psi_f$. 
Here we used Lagrange multipliers $\mu_m$ and $\mu_f$ to impose the constraints for multiscale basis construction. 
We remark that we have used the notations $\psi_m,\psi_f,\mu_m,\mu_f$ to denote the vector representations of the corresponding functions in terms of fine scale basis. 
For example $\psi_m$ is the vector of coefficients of the matrix pressure expanded in terms of fine scale basis. 
We set $F_m = \delta_{i,j}$ and $F_f = 0$ for construction of multiscale basis function for porous matrix $\psi^0 = (\psi^0_m, \psi^0_f)$. For multiscale basis function for fracture network, we set  $F_m =  0$ and $F_f = \delta_{i,j}\delta_{m,l}$. In Figure \ref{fig:msbf20}, we depict multiscale basis functions for oversampled region $K^+_i = K^4_i$ (four oversampling coarse cell layers) on coarse mesh $20 \times 20$.

\textbf{Multiscale basis function for displacements.} 
The construction is similar to that of pressure. More precisely, 
we construct a set of basis functions $\psi^{X,i} :=  (\psi^{X,i}_{x}, \psi^{X,i}_{y})$ 
and $\psi^{Y,i} := (\psi^{Y,i}_{x}, \psi^{Y,i}_{y})$, which minimize the energy for elasticity problem operator  restricted in the region $K_i^+$
and satisfy the constraints described below for all $K_j \subset K_i^+$:\\
(1)  X-component, $\psi^{X,i}$ :
\[
\int_{K_j} \psi_{x}^{X,i} \, dx = \delta_{i,j}, \quad
\int_{K_j} \psi_{y}^{X,i} \, dx = 0,
\]
(2)  Y-component, $\psi^{Y,i}$ :
\[
\int_{K_j} \psi_{x}^{Y,i} \, dx = 0, \quad 
\int_{K_j} \psi_{y}^{Y,i} \, dx = \delta_{i,j}.
\]
For further error reduction, we can add additional basis function for heterogeneous source term.  

{
This results in solving the following local problems in $K_i^+$:
\begin{equation}
\label{eq:basis-u}
\begin{pmatrix}
D_x^{K^+_i} & D_{xy}^{K^+_i} & S^T_x & 0 \\
D_{xy}^{K^+_i} & D_{yy}^{K^+_i} & 0 &  S^T_y \\
S_x & 0 & 0 & 0 \\
0 & S_y & 0 & 0 \\
\end{pmatrix} 
\begin{pmatrix}
\psi_x \\
\psi_y \\
\mu_x \\
\mu_y \\
\end{pmatrix} = 
\begin{pmatrix}
0 \\
0 \\
F_x \\
F_y \\
\end{pmatrix}
\end{equation}
with  zero Dirichlet boundary conditions on $\partial K^+_i$ for $\psi_x$ and $\psi_y$. 
We set $(F_x, F_y) = (\delta_{i,j}, 0)$ and $(0, \delta_{i,j})$ for construction of multiscale basis function for X and Y-components. 
}
In Figure \ref{fig:msbf20}, we depict multiscale basis functions for displacements in oversampled domain $K^+_i = K^4_i$. 

{
In general, the permeability or elastic coefficients can be heterogeneous, where for high-construct cases more basis should be used and constrained energy minimization (CEM) GMsFEM can identify important modes  \cite{chung2017constraint}. 
}

{
We note that, the fracture contributions are divided in each coarse cell and then coupled. Each local fracture network introduce an additional degree of freedom for current coarse cell. In general, CEM-GMsFEM can be applied, where local spectral problem automatically identify important modes \cite{chung2017constraint}. 
}

\textbf{Coarse scale coupled system. }
We first define a projection matrix using the multiscale basis functions
\[
R = \begin{pmatrix}
R_{mm} & R_{mf} & 0  & 0 \\
R_{fm} & R_{ff} & 0  & 0 \\
0  & 0 & R_{xx} & R_{xy} \\
0  & 0 & R_{yx} & R_{yy}  \\
\end{pmatrix}, 
\]
where
\[
R_{mm}^T = \left[ 
\psi^{0,0}_m, \psi^{1,0}_m  \ldots \psi^{N_c,0}_m
\right],\quad
 R_{ff}^T = \left[ 
\psi^{0,1}_f \ldots \psi^{0,L_0}_f, 
\psi^{1,1}_f \ldots \psi^{1,L_1}_f,
\ldots, 
\psi^{N_c,1}_f  \ldots \psi^{N_c,L_{N_c}}_f
 \right],
\]\[
R_{mf}^T = \left[ 
\psi^{0,0}_f, \psi^{1,0}_f \ldots \psi^{N_c,0}_f
 \right], \quad
R_{fm}^T = \left[ 
\psi^{0,1}_m \ldots \psi^{0,L_0}_m, 
\psi^{1,1}_m \ldots \psi^{1,L_1}_m,
\ldots, 
\psi^{N_c,1}_m  \ldots \psi^{N_c,L_{N_c}}_m
 \right], 
\]\[
R_{xx}^T = \left[ \psi^{X,0}_{x}, \ldots \psi^{X,N_c}_{x} \right]
 \quad
R_{xy}^T = \left[ \psi^{X,0}_{y}, \ldots \psi^{X,N_c}_{y} \right]
\quad
R_{yx}^T = \left[ \psi^{Y,0}_{x},  \ldots \psi^{Y,N_c}_{x} \right]
\quad
R_{yy}^T = \left[ \psi^{Y,0}_{y},  \ldots  \psi^{Y,N_c}_{y} \right].
\]
In the above definition, $\psi_m^{i,l}$ is the basis function for matrix pressure corresponding to the coarse block $K_i$ and the continuum $l$. The
definition for $\psi_f^{i,l}$ is the basis function for fracture pressure corresponding to the coarse block $K_i$ and the continuum $l$.
The notation $\psi_m^{i,l}$ stands for both the function and its vector representation in fine grid basis. 
{We note that we construct only decoupled multiscale basis functions for flow and mechanics. Coupled construction of the multiscale basis functions can provide better results and will be considered  and investigated in the future works.}

Finally, we obtain following upscaled coarse grid model
\begin{equation}
\label{t-nlmc2}
\left(\frac{1}{\tau} \bar{M} + \bar{A} \right) \bar{y} = \bar{F},
\end{equation}
where 
$\bar{A} = R A R^T$, $\bar{F} = R F$, 
$\bar{y} = (\bar{p}_m, \bar{p}_f, \bar{u}_x, \bar{u}_y)$. 
Here
$\bar{p}_m$, $\bar{p}_f$, $\bar{u}_x$, $\bar{u}_y$ are the average  solution on coarse grid cell for matrix, fracture, displacement X and Y components, respectively.
For mass matrix, we can use a property of the constructed multiscale basis functions, and obtain diagonal mass matrix by direct calculation on the coarse grid 
\[
\bar{M} = 
\begin{pmatrix}
\bar{M}_m & 0 & 0 & 0 \\
0 & \bar{M}_f & 0 & 0\\
0 & 0 & 0 & 0\\
0 & 0 & 0 & 0\\
\end{pmatrix}, 
\] 
where 
$\bar{M}_m = \text{diag}\{ a_m |K_i| \}$, 
$\bar{M}_f = \text{diag}\{ a_f |\gamma_i| \}$. 
The coarse grid upscaled model has only one coarse degree of freedom (DOF) for each fracture network and  provides an effective coarse scale model with physical meaning, and leads to a fast and accurate solver for the coupled poroelasticity problem.

\section{Numerical results}

\begin{figure}[h!]
\centering
\includegraphics[width=0.9 \textwidth]{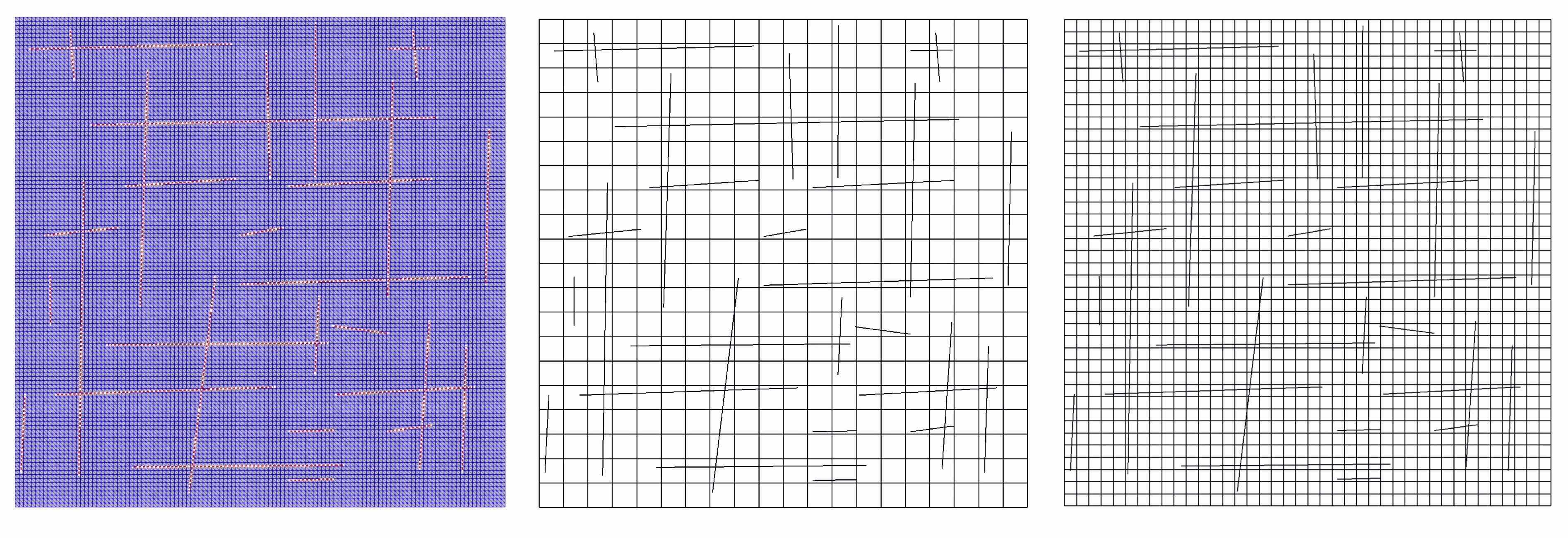}
\caption{Computational grids  with 30 fracture lines.  
First: Coarse grid $20 \times 20$ with 400 cells.
Second: Coarse grid $40 \times 40$ with 1600 cells.
Third: Fine grid for matrix domain $\Omega$ with 14641 vertices and 28800 cells (blue). Fine gird for fracture domain $\gamma$ with 1042 cells (red and white)}
\label{fig:mesh3}
\end{figure}

\begin{figure}[h!]
\centering
\includegraphics[width=0.6 \textwidth]{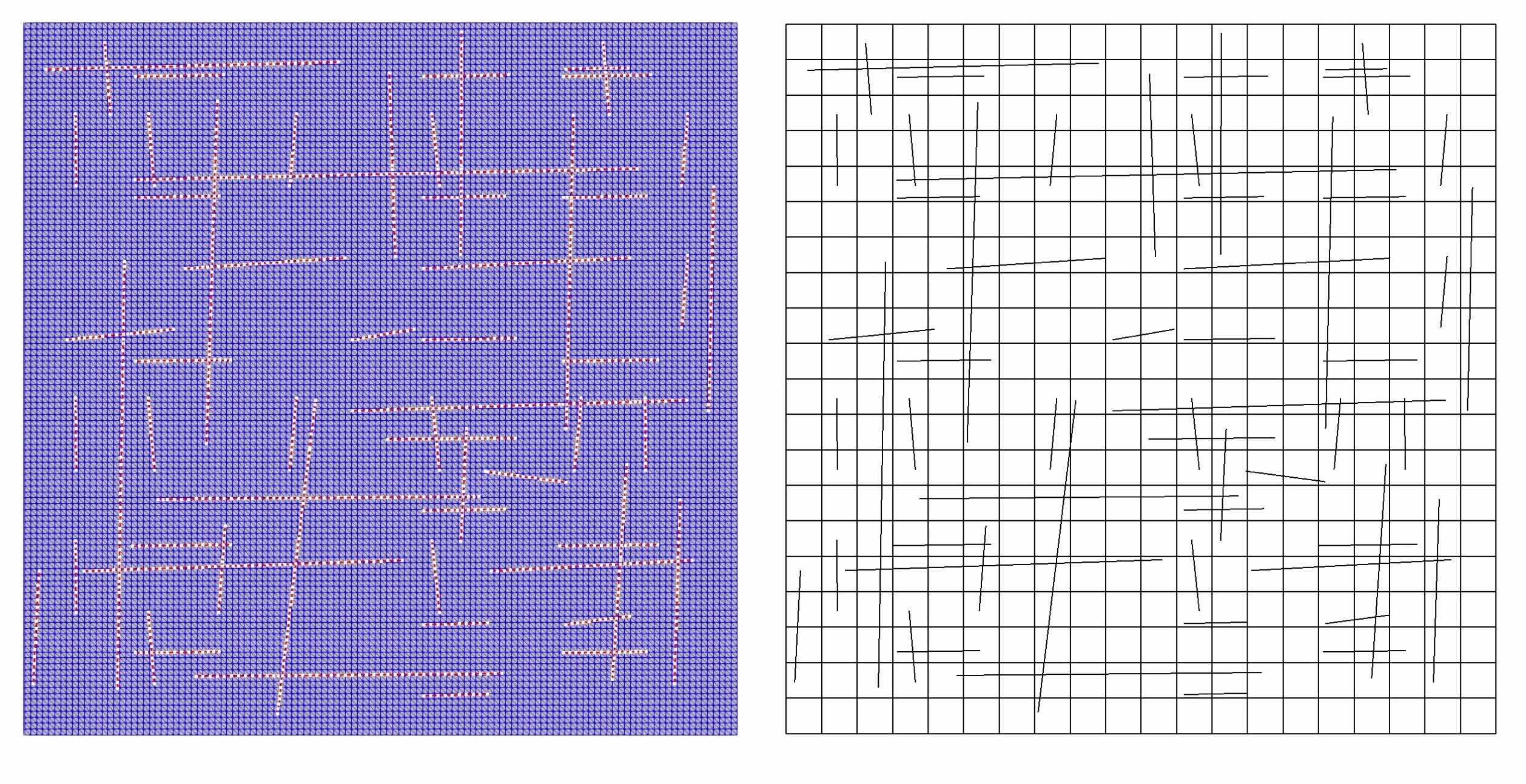}
\caption{Computational grids with 60 fracture lines.  
First: Coarse grid $20 \times 20$ with 400 cells.
Second: Fine grid for matrix domain $\Omega$ with 14641 vertices and 28800 cells (blue). Fine gird for fracture domain $\gamma$ with 1312 cells (red and white)}
\label{fig:mesh3b}
\end{figure}

We present  numerical results  for poroelastic model in $\Omega$ with length of 1 meter in both directions. 
We consider two test cases: (1) domain with 30 fractures and (2) domain with 60 fractures. In Figures \ref{fig:mesh3} and \ref{fig:mesh3b}, we show computational coarse and fine grids, where the fractures are depicted with red color and fine mesh with blue color. 
For fracture  network, we constructed separate mesh and for domain $\Omega$, we use structured fine mesh.  
We consider two coarse grids with 400 cells and with 1600 cells.
The coarse grids are uniform.

For coupled poroelastic model, we use following parameters:
\begin{itemize}
\item Elastic parameters: $\mu = \frac{E}{2 (1 + \nu)}$ and  $\lambda = \frac{E \nu}{(1+ \nu) ( 1- 2 \nu)}$, where $E = 10 \times 10^9$,  $\nu = 0.3$ and $\alpha = 0.1$,
\item Flow parameters $a_m = 10^{-6}$, $a_f = 10^{-7}$, {$b_m = 10^{-11}$, $b_f = 10^{-6}$ and $\beta = 10^{-10}$.}
\end{itemize}
Boundary condition for the displacement: 
$u_x = 0.0$ on the left and right boundaries, 
$u_y = 0.0$ on the bottom and top.
We set a point source at the two coarse cells with $q = 0.01$ and set initial pressure  $p_0 = 10^7$. 
We simulate $t_{max} = 10$ years with 50 time steps for multiscale and fine scale solvers. 

\begin{figure}[h!]
\centering
\includegraphics[width=0.75 \textwidth]{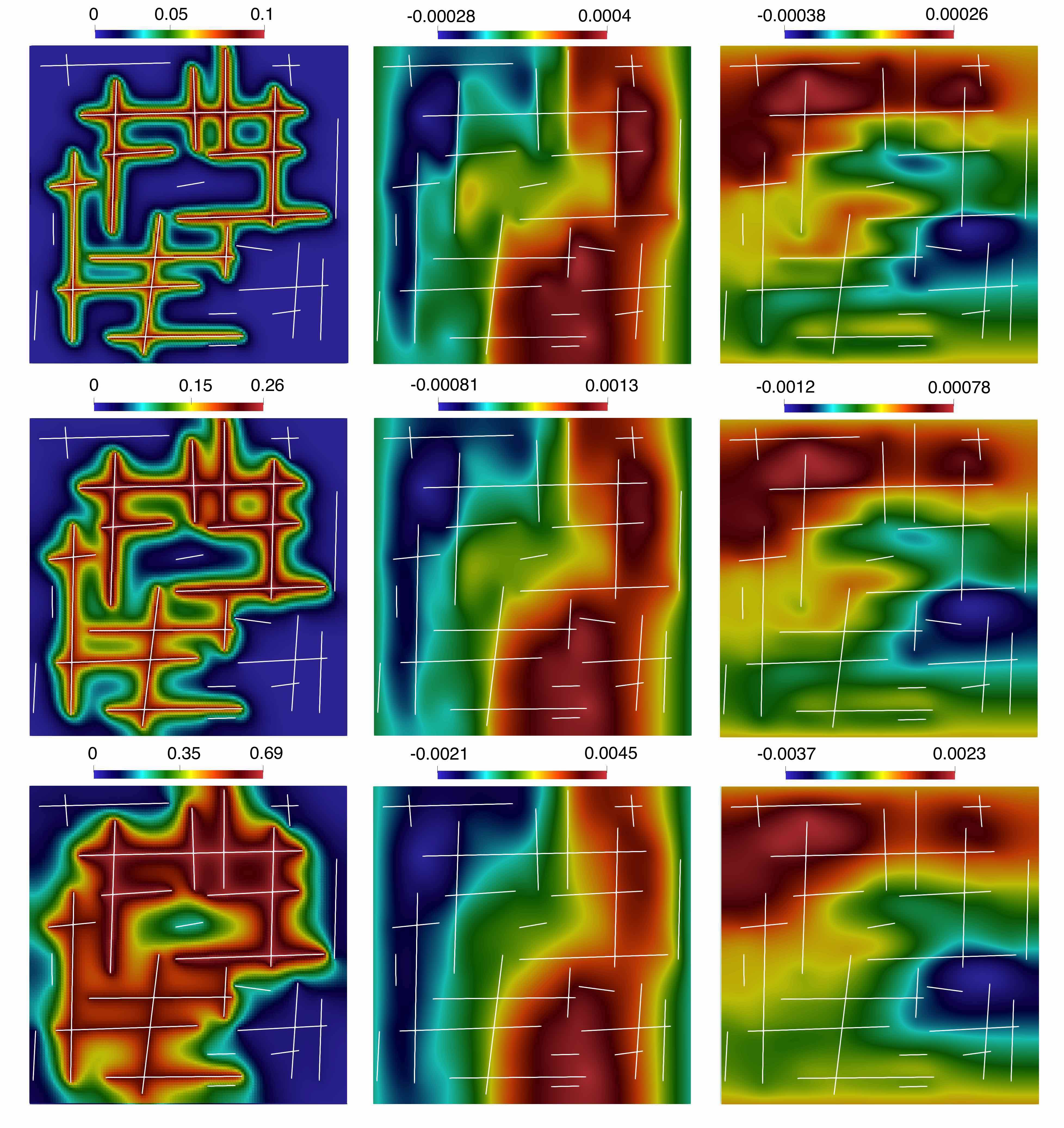}
\caption{Fine scale solution for pressure ($p^* = (p - p_0)/p_0$) and displacements (from left to right) for the different time layers $t_5$, $t_{15}$ and $t_{50}$ (from top to bottom). Test case with 30 fractures.
}
\label{fig:s3}
\end{figure}

\begin{figure}[h!]
\centering
\includegraphics[width=1.0 \textwidth]{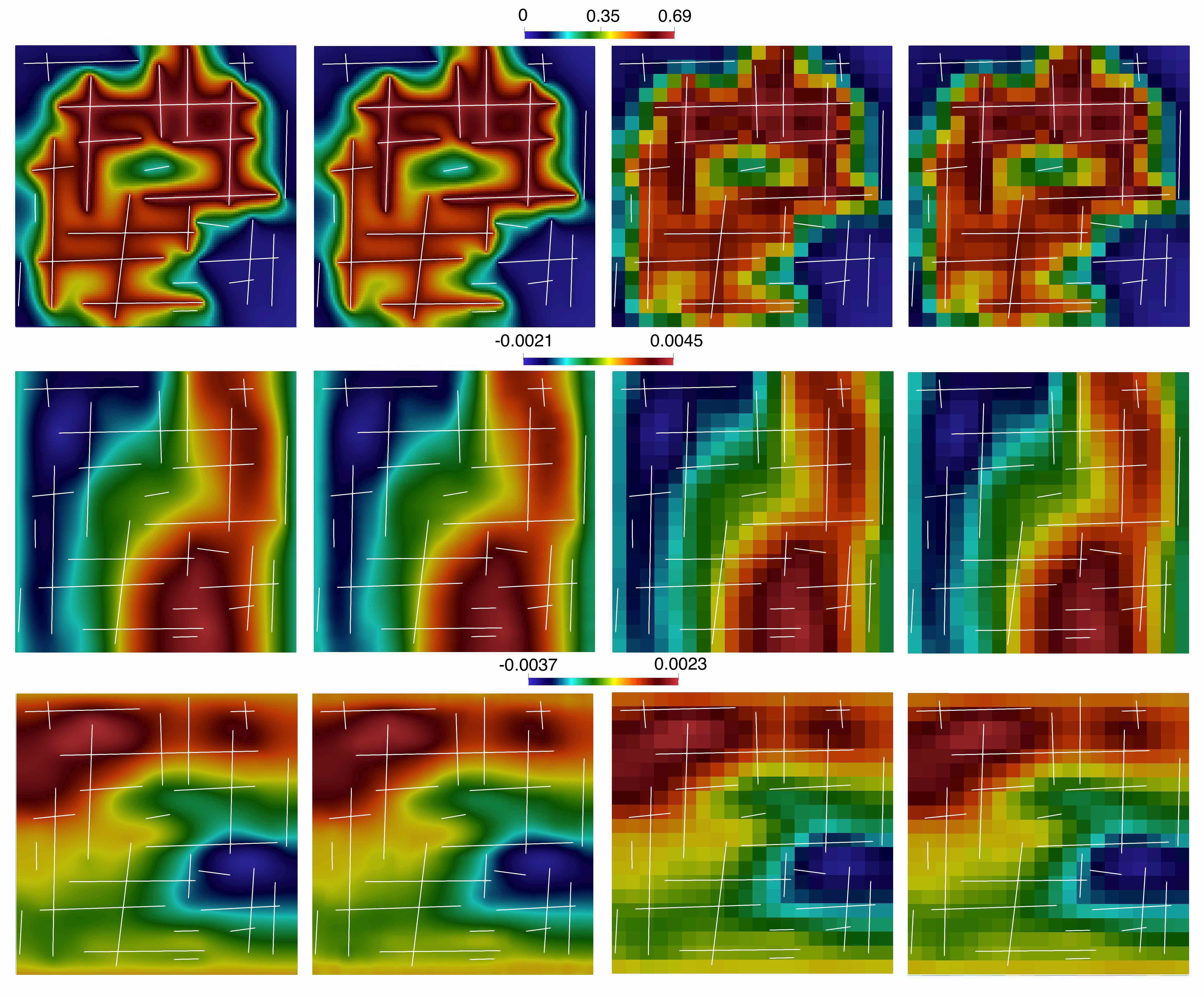}
\caption{Numerical results for pressure   ($p^* = (p - p_0)/p_0$) and displacements at final time.  Test case with 30 fractures.
First row: pressure, $p_m$. 
Second row: displacement, $u_x$. 
Third row: displacement, $u_y$. 
First column: fine scale solution with $DOF_f = 59124$. 
Second column: reconstructed fine scale solution from upscaled coarse grid solution with $DOF_c = 1393$. 
Third column: coarse cell average for fine scale solution.  
Fourth column: coarse cell average for upscaled coarse grid. Coarse grid $20 \times 20$ ($K^+ = K^4$).}
\label{fig:ms3}
\end{figure}

\begin{figure}[h!]
\centering
\includegraphics[width=0.75 \textwidth]{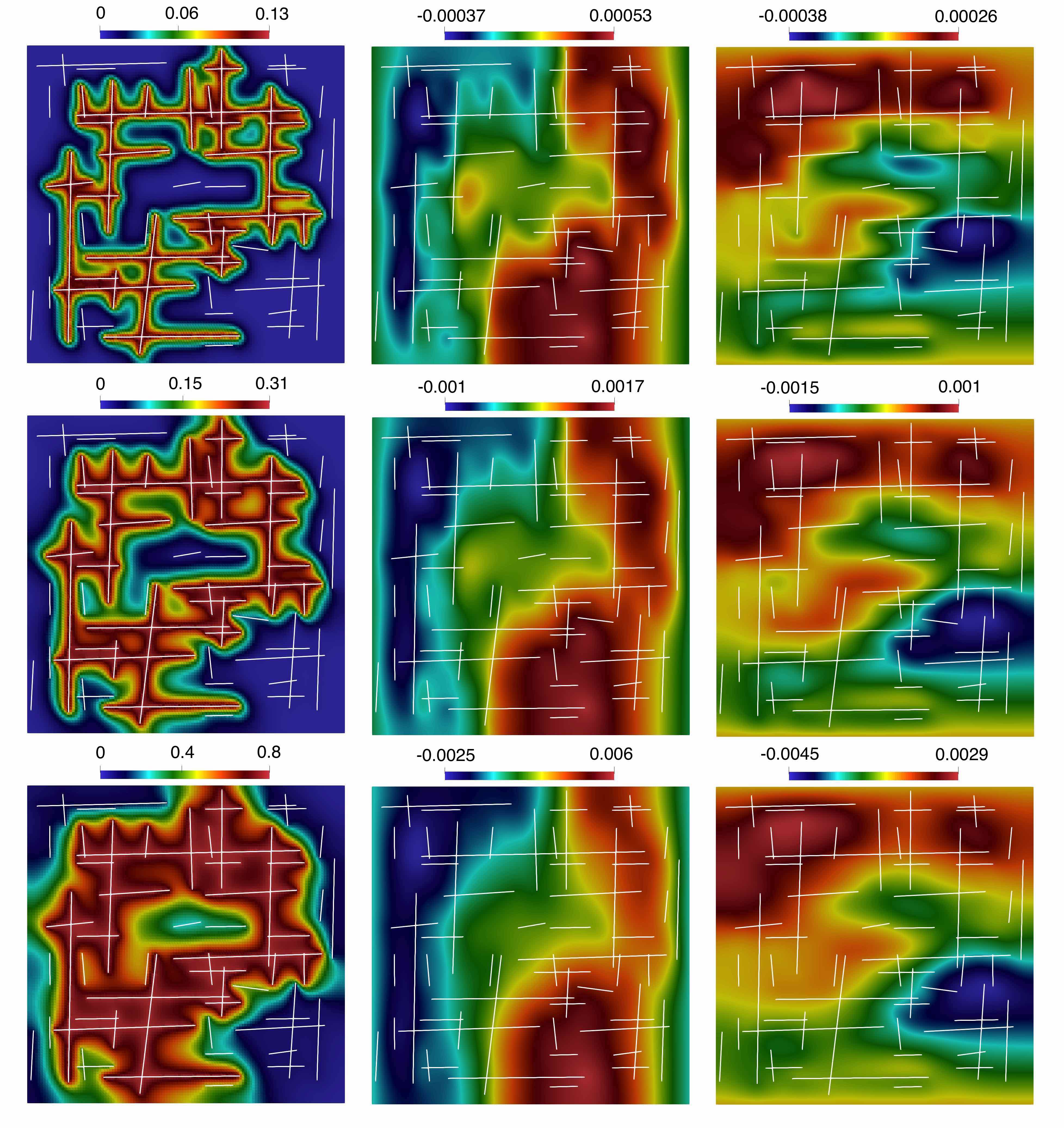}
\caption{Fine scale solution for pressure ($p^* = (p - p_0)/p_0$) and displacements (from left to right) for the different time layers $t_5$, $t_{15}$ and $t_{50}$ (from top to bottom). Test case with 60 fractures.
}
\label{fig:s3b}
\end{figure}

\begin{figure}[h!]
\centering
\includegraphics[width=1 \textwidth]{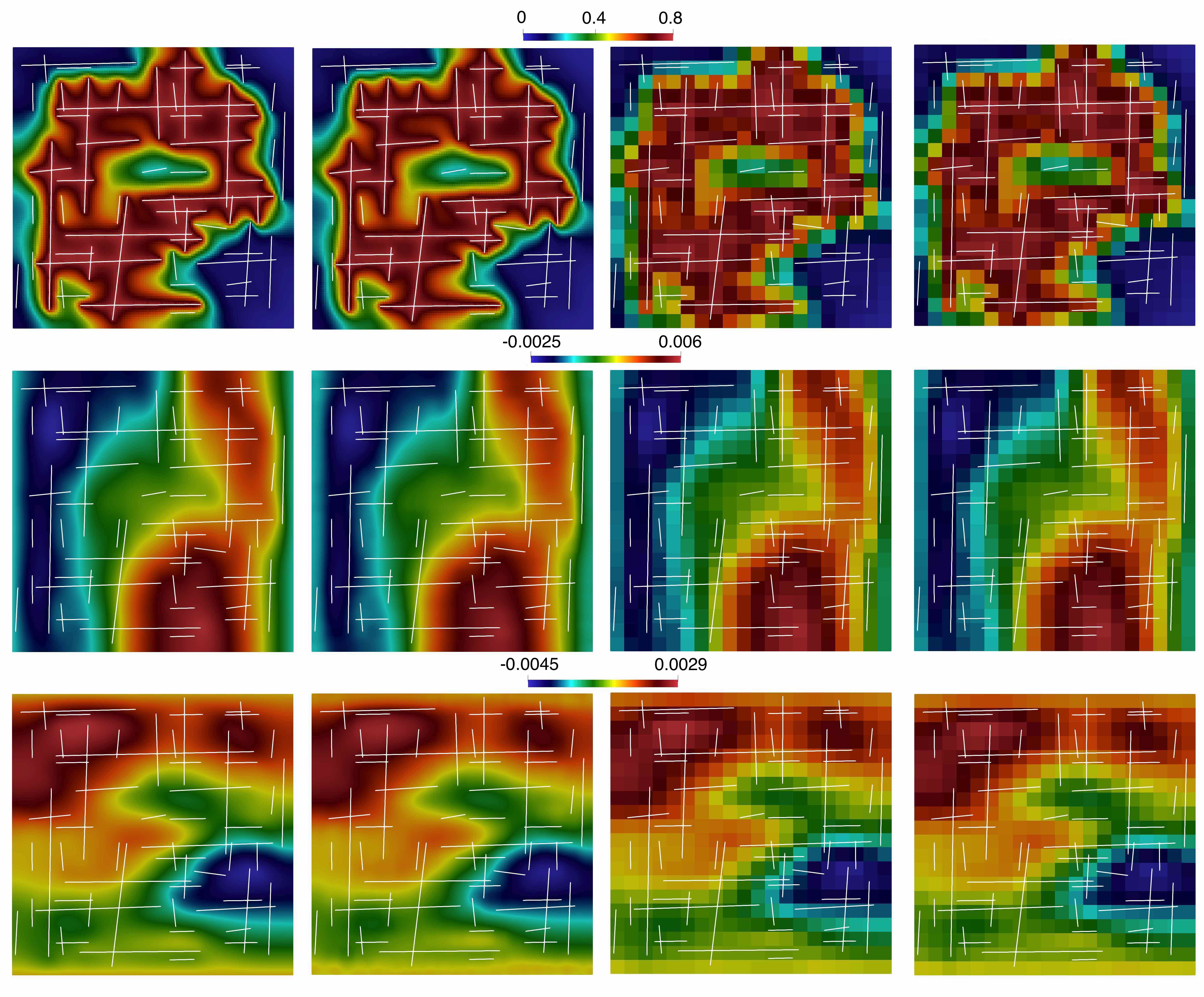}
\caption{Numerical results for pressure  ($p^* = (p - p_0)/p_0$) and displacements at final time.  Test case with 60 fractures.
First row: pressure, $p_m$. 
Second row: displacement, $u_x$. 
Third row: displacement, $u_y$. 
First column: fine scale solution with $DOF_f = 59394$. 
Second column: reconstructed fine scale solution from upscaled coarse grid solution with $DOF_c = 1484$. 
Third column: coarse cell average for fine scale solution.  
Fourth column: coarse cell average for upscaled coarse grid. Coarse grid $20 \times 20$ ($K^+ = K^4$)}
\label{fig:ms3b}
\end{figure}

\begin{table}[h!]
\begin{center}
\begin{tabular}{|c|ccc|}
\hline
$K^s$ & $e_p$ & $e_{u_x}$ & $e_{u_y}$ \\
\hline
\multicolumn{4}{|c|}{Coarse grid $20 \times 20$}  \\
\hline
1		&	4.740	&	 86.865	&	 82.598 \\
2		&	0.723	&	 43.721	&	 37.034 \\
3		&	0.369	&	 6.716		&	 4.668 \\
4		&	0.359	&	 2.718		&	 2.854 \\
\hline
\end{tabular}
\,\,\,\,
\begin{tabular}{|c|ccc|}
\hline
$K^s$ & $e_p$ & $e_{u_x}$ & $e_{u_y}$ \\
\hline
\multicolumn{4}{|c|}{Coarse grid $40 \times 40$}  \\
\hline
1		&	1.986	&	 96.667	&	 95.454 \\
2		&	0.191	&	 78.718	&	 74.957 \\
3		&	0.174	&	 30.550	&	 25.220 \\
4		&	0.158	&	 4.1302	&	 3.321 \\
6		&	0.157	&	 1.127		&	 1.233 \\
\hline
\end{tabular}
\end{center}
\caption{Numerical results of relative errors (\%) at the final simulation time. $DOF_f= 59124$ and $DOF_c = 1393$. Test case with 30 fractures.}
\label{tab:errs3}
\end{table}

\begin{figure}[h!]
\centering
\includegraphics[width=0.49 \textwidth]{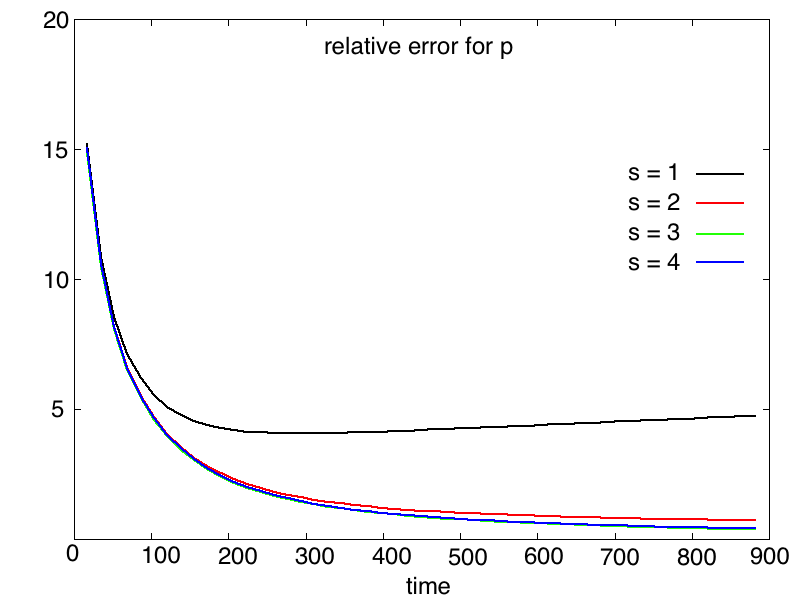}
\includegraphics[width=0.49 \textwidth]{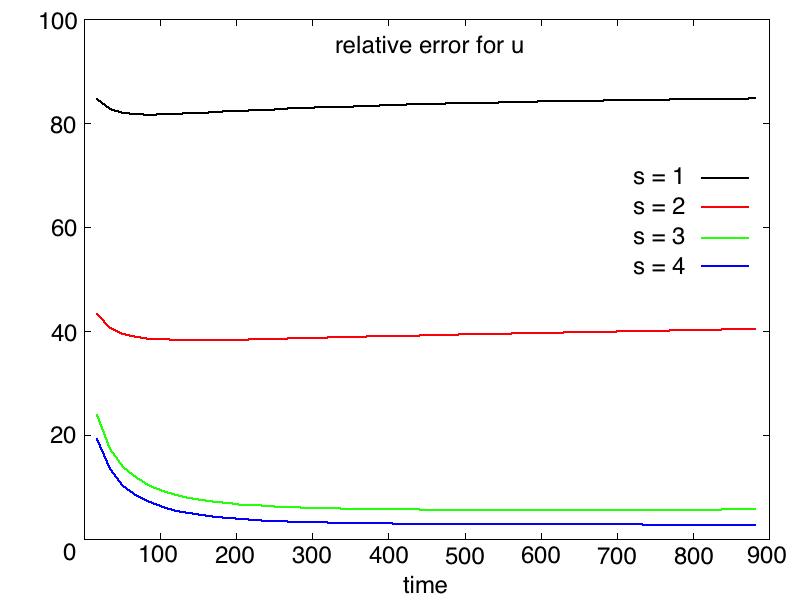}
\caption{Relative errors by time for coarse mesh $20 \times 20$ with $K^s$. Test case with 30 fractures.}
\label{fig:err-t3-20}
\end{figure}


\begin{figure}[h!]
\centering
\includegraphics[width=0.49 \textwidth]{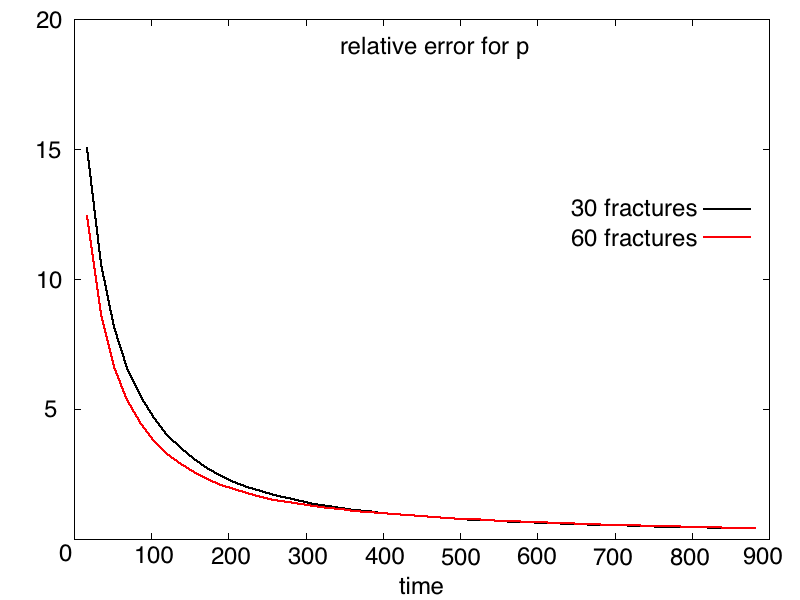}
\includegraphics[width=0.49 \textwidth]{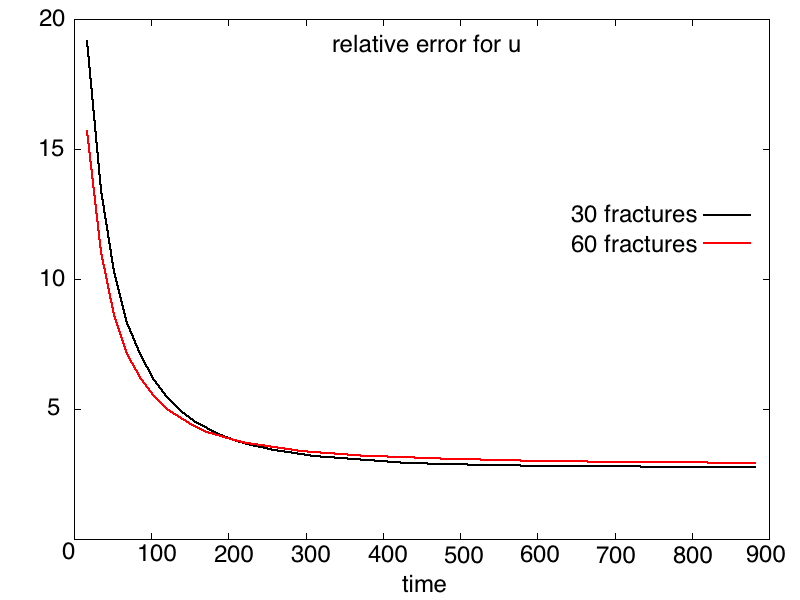}
\caption{Relative errors by time for coarse mesh $20 \times 20$ with $K^4$. Test cases with 30 and 60 fractures.}
\label{fig:err-t3b}
\end{figure}

We use $DOF_c$ to denote problem size of the coarse-grid upscaled model and $DOF_f$ for the fine grid system size.   
{
To compare the results, we use the relative $L^2$ error between coarse cell average of the fine-scale solution 
$\bar{p}^{fine}_m, \bar{u}^{fine}_{x}, \bar{u}^{fine}_{y}$ 
and upscaled coarse grid  solutions 
$\bar{p}_m, \bar{u}_x, \bar{u}_y$
\begin{equation}
\begin{split}
e_p = ||\bar{p}^{fine}_{m} - \bar{p}_m ||_{L^2}, \quad 
e_{u_x} = ||\bar{u}_{x}^{fine} - \bar{u}_x ||_{L^2}, \quad 
e_{u_y} = ||\bar{u}_{y}^{fine} - \bar{u}_y ||_{L^2}, \\
|| \bar{v}_f - \bar{v} ||^2_{L^2} =  
\frac{ \sum_K (\bar{v}^K_f - \bar{v}^K)^2}{\sum_K (\bar{v}^K_f)^2 }, \quad
\bar{v}^K_f = \frac{1}{|K|} \int_K v_f \, dx, \quad v = p, u_x, u_y,
\end{split}
\end{equation}
for matrix pressure and displacements.
}

Fine grid solution for computational domain with 30 fractures is presented in Figure \ref{fig:s3} for the different time instants $t_5$, $t_{15}$ and $t_{50}$, where $t_n = n\tau$.
On the first column of the figure, we depict pressure $p^* = (p - p_0)/p^0$, on the second and third row columns -- displacements $u
_x$ and $u_y$.
Comparison of the fine grid and coarse grid upscaled solutions are presented in Figure \ref{fig:ms3} at final time. We perform computations on the coarse grid with 400 cells with $4$ oversampling layers in the construction of basis functions ($K^+ = K^4$). In the first column, we depict fine grid pressure solution; in the second column -- reconstructed fine scale solution from upscaled coarse grid solution, in the third column -- coarse cell average for fine scale solution and in the fourth column -- coarse cell average for upscaled coarse grid. Fine grid system has size $DOF_f = 59124$. By performing NLMC method, we reduce size of system to $DOF_c = 1393$.
At final time, we have less than one percent of error for pressure and near $2.5 \%$ for displacement.

In Table \ref{tab:errs3}, we present relative errors at final time for two coarse grids and for different numbers of oversampling layers for the oversample region $K^s$ with $s = 1,2,3,4$ and $6$,
where $K^s$ is obtained by extending $K$ by $s$ coarse grid layers. 
From the numerical results, we observe a good convergence behavior, when we take sufficient number of oversampled layers. For the coarse mesh with 400 cells, when we take 4 oversampling layers, we have $0.359 \%$ relative error for pressure, for displacement -- $2.718 \%$ ($u_x$) and $2.854 \%$ ($u_y$).
For the coarse mesh with 1600 cells with 6 oversampling layers, relative error is $0.157 \%$ for pressure, for displacement -- $1.127 \%$ ($u_x$) and $1.233 \%$ ($u_y$).
We note that, on the $20 \times 20$ coarse mesh, the size of upscaled system is $DOF_c = 1393$ and for the $40 \times 40$ coarse mesh, we have $DOF_c = 5165$. 
{ From the Table \ref{tab:errs3}, we observe that we can use smaller number of oversampling layers for pressure than for displacements. For pressure is enough to take 2 oversampling layers for obtaining errors smaller than one percent for both coarse grids, on the other hand for displacements we should take 4 or 6 oversampling layers in coarse grids $20 \times 20$ and $40 \times 40$, respectively. Note tat, in general the presented algorithm can work with different numbers of oversampling layers for pressure and displacement due to coupled construction of the coarse grid system.}

In Figure \ref{fig:err-t3-20}, we present relative errors for pressure and displacements vs time with different number of oversampling layers $K^s$, $s = 1,2,3,4, 6$.
All results show good accuracy of the proposed method for coupled poroelasticity problems in fractured media.

Next, we consider test case with 60 fractures.
In Figure \ref{fig:s3b}, we shown solution of the problem for the different time layers $t_5$, $t_{15}$ and $t_{50}$.
Comparison of the fine grid and coarse grid upscaled solutions are presented in Figure \ref{fig:ms3b} at final time for coarse grid with 400 cells and $4$ oversampling layers ($K^+ = K^4$). In the first column, we depict fine grid pressure solution; in the second column -- reconstructed fine scale solution from upscaled coarse grid solution, in the third column -- coarse cell average for fine scale solution and in the fourth column -- coarse cell average for upscaled coarse grid.
Fine grid system has size $DOF_f = 59394$. By performing NLMC method, we reduce size of system to $DOF_c = 1484$.
At final time, we have $0.4217 \%$ of error for pressure and for displacement -- $2.739 \%$ ($u_x$) and $3.124 \%$ ($u_y$).
In Figure \ref{fig:err-t3b}, we depict relative errors vs time for $K^4$ on coarse mesh $20 \times 20$ for test case with 30 and 60 fractures.
We obtain similar results with good accuracy for both test cases.

Next, we discuss the computational advantages of our approach. The computatonal time is divided into offline and online stages. In offline stage (preprocessing), we generate local domains, calculate multiscale basis functions and generate coarse grid system. In online stage, we solve coarse grid problem, with different imput parameters (source term, boundary conditions, time steps, etc.).
Let $DOF_f$ is the size of fine scale solution $y = (p_m, p_f, u)$, then the dimension of the fine grid coupled problem is $DOF_f \times DOF_f$.
The coarse grid system size is $DOF_c$ for coupled poroelasticity problem, that depends on the coarse grid size and the number of local multiscale basis functions. In each local domain (coarse grid cell), we have degree of freedom for displacement X and Y components, vof matrix pressure and additional degree of freedom for each fracture network in current coarse cell. We note that, the number of degree of freedom is similar to classic embedded fracture model (EFM).
For two dimensional problems ,we have $M = 3 + M_i$ ($(p_m, p^1_f... p_f^{M_i}, u_x, u_y)$) degrees of freedoms in local domain, where $M_i$ is the number of the fracture networks in coarse cell $K_i$. Therefore, the size of coarse grid system is $DOF_c = \sum_{K_i}(3 + M_i)$.

Let $N_f$ is the number of cells for fracture network mesh, $N_c$ and $N_v$ are the number of cells and vertices on fine grid for domain $\Omega$. Then for two dimensional problems with finite element approximation for displacements equation and finite volume approximation for flow problem, we have $DOF_f = N_f + N_c + 2 N_v$.
Then, we can compare a computational cost of solving coarse and fine grid problems.
For example in test case with 30 fractures, a coarse solution has $DOF_c = 1393$ in the coarse grid with $400$ cells, where we have $400$ and $800$ degrees of freedom for matrix pressure and displacement X and Y components, and $193$ degrees of freedom for fractures.
For the fine scale system $DOF_f = 59124$ on fine grid with $N_v = 14641$ vertices and $N_c = 28800$ cells.
Then, we can obtain accurate solution for multiscale solver using only $2.3 \%$ from $DOF_f$. We note that, the number of $M_i$ in $K_i$ and therefore size of coarse grid system is independent on fine grid size and  a few basis functions can approximate the fine scale solution
accurately no matter how fine is the fine grid. When we use classic direct solver, the solution time of the time dependent coupled fine grid problem is $81.17$ seconds and $5.74$ seconds for coarse grid.
We have computational gain in the simulations, because in each time step, the proposed method solves a small coarse grid system compared to the fine-grid system.

\section{Conclusion}
In this paper, our goal is to develop an upscaled model for a poroelastic system 
in fractured media. There are several contributions.
First, we construct an embedded fracture model for a coupled flow and mechanics system. Secondly, 
based on this system, we develop
a nonlocal upscaled model for efficient numerical simulations. The construction of the upscaled 
model is motivated by the NLMC method.
The main idea is to construct basis functions for each continuum within a 
local coarse region such that the resulting coarse degrees of freedom have physical meanings.
Moreover, these basis functions have decay property
thanks to an energy minimization principle, which can guarantee an accurate 
approximation of the solution. We have presented several numerical 
tests to show that our upscaled model can give
accurate solutions with a small computational cost. 

\section*{Acknowledgements}
MV's  work is supported by the grant of the Russian Scientific Found N17-71-20055.  YE's  is supported by the mega-grant of the Russian Federation Government (N 14.Y26.31.0013).  
EC's work is partially supported by Hong Kong RGC General Research Fund (Project 14304217)
and CUHK Direct Grant for Research 2017-18.

\bibliographystyle{plain}
\bibliography{poroelastic}

\begin{thebibliography}{10}

\bibitem{akkutlu2018multiscale}
I~Yucel Akkutlu, Yalchin Efendiev, Maria Vasilyeva, and Yuhe Wang.
\newblock Multiscale model reduction for shale gas transport in poroelastic
  fractured media.
\newblock {\em Journal of Computational Physics}, 353:356--376, 2018.

\bibitem{akkutlu2015multiscale}
IY~Akkutlu, Yalchin Efendiev, and Maria Vasilyeva.
\newblock Multiscale model reduction for shale gas transport in fractured
  media.
\newblock {\em Computational Geosciences}, pages 1--21, 2015.

\bibitem{barenblatt1960basic}
GI~Barenblatt, Iu~P Zheltov, and IN~Kochina.
\newblock Basic concepts in the theory of seepage of homogeneous liquids in
  fissured rocks [strata].
\newblock {\em Journal of applied mathematics and mechanics}, 24(5):1286--1303,
  1960.

\bibitem{bosma2017multiscale}
Sebastian Bosma, Hadi Hajibeygi, Matei Tene, and Hamdi~A Tchelepi.
\newblock Multiscale finite volume method for discrete fracture modeling on
  unstructured grids (ms-dfm).
\newblock {\em Journal of Computational Physics}, 2017.

\bibitem{brown2016generalized}
Donald~L Brown and Maria Vasilyeva.
\newblock A generalized multiscale finite element method for poroelasticity
  problems i: linear problems.
\newblock {\em Journal of Computational and Applied Mathematics}, 294:372--388,
  2016.

\bibitem{brown2016generalized2}
Donald~L Brown and Maria Vasilyeva.
\newblock A generalized multiscale finite element method for poroelasticity
  problems ii: Nonlinear coupling.
\newblock {\em Journal of Computational and Applied Mathematics}, 297:132--146,
  2016.

\bibitem{CELV2015}
E.~T. Chung, Y.~Efendiev, G.~Li, and M.~Vasilyeva.
\newblock Generalized multiscale finite element method for problems in
  perforated heterogeneous domains.
\newblock {\em to appear in Applicable Analysis}, 255:1--15, 2015.

\bibitem{chung2016adaptive}
Eric Chung, Yalchin Efendiev, and Thomas~Y Hou.
\newblock Adaptive multiscale model reduction with generalized multiscale
  finite element methods.
\newblock {\em Journal of Computational Physics}, 320:69--95, 2016.

\bibitem{chung2017coupling}
Eric~T Chung, Yalchin Efendiev, Tat Leung, and Maria Vasilyeva.
\newblock Coupling of multiscale and multi-continuum approaches.
\newblock {\em GEM-International Journal on Geomathematics}, 8(1):9--41, 2017.

\bibitem{chung2017constraint}
Eric~T Chung, Yalchin Efendiev, and Wing~Tat Leung.
\newblock Constraint energy minimizing generalized multiscale finite element
  method.
\newblock {\em Computer Methods in Applied Mechanics and Engineering},
  339:298--319, 2018.

\bibitem{chung2017non}
Eric~T Chung, Yalchin Efendiev, Wing~Tat Leung, Yating Wang, and Maria
  Vasilyeva.
\newblock Non-local multi-continua upscaling for flows in heterogeneous
  fractured media.
\newblock {\em arXiv preprint arXiv:1708.08379}, 2017.

\bibitem{coussy2004poromechanics}
Olivier Coussy.
\newblock {\em Poromechanics}.
\newblock John Wiley \& Sons, 2004.

\bibitem{Quarteroni2008coupling}
Carlo D'angelo and Alfio Quarteroni.
\newblock On the coupling of 1d and 3d diffusion-reaction equations:
  application to tissue perfusion problems.
\newblock {\em Mathematical Models and Methods in Applied Sciences},
  18(08):1481--1504, 2008.

\bibitem{douglas1990dual}
Jim Douglas~Jr and T~Arbogast.
\newblock Dual porosity models for flow in naturally fractured reservoirs.
\newblock {\em Dynamics of Fluids in Hierarchical Porous Media}, pages
  177--221, 1990.

\bibitem{d2012mixed}
Carlo D’Angelo and Anna Scotti.
\newblock A mixed finite element method for darcy flow in fractured porous
  media with non-matching grids.
\newblock {\em ESAIM: Mathematical Modelling and Numerical Analysis},
  46(2):465--489, 2012.

\bibitem{EGG_MultiscaleMOR}
Y.~Efendiev, J.~Galvis, and E.~Gildin.
\newblock Local-global multiscale model reduction for flows in highly
  heterogeneous media.
\newblock {\em Journal of Computational Physivs}, 231 (24):8100--8113, 2012.

\bibitem{egh12}
Y.~Efendiev, J.~Galvis, and T.~Hou.
\newblock Generalized multiscale finite element methods.
\newblock {\em Journal of Computational Physics}, 251:116--135, 2013.

\bibitem{eh09}
Y.~Efendiev and T.~Hou.
\newblock {\em {Multiscale Finite Element Methods: Theory and Applications}},
  volume~4 of {\em Surveys and Tutorials in the Applied Mathematical Sciences}.
\newblock Springer, New York, 2009.

\bibitem{efendiev2015hierarchical}
Yalchin Efendiev, Seong Lee, Guanglian Li, Jun Yao, and Na~Zhang.
\newblock Hierarchical multiscale modeling for flows in fractured media using
  generalized multiscale finite element method.
\newblock {\em arXiv preprint arXiv:1502.03828}, 2015.
\newblock to appear in International Journal on Geomathematics, (DOI)
  10.1007/s13137-015-0075-7.

\bibitem{formaggia2014reduced}
Luca Formaggia, Alessio Fumagalli, Anna Scotti, and Paolo Ruffo.
\newblock A reduced model for darcy’s problem in networks of fractures.
\newblock {\em ESAIM: Mathematical Modelling and Numerical Analysis},
  48(4):1089--1116, 2014.

\bibitem{girault2016convergence}
Vivette Girault, Kundan Kumar, and Mary~F Wheeler.
\newblock Convergence of iterative coupling of geomechanics with flow in a
  fractured poroelastic medium.
\newblock {\em Computational Geosciences}, 20(5):997--1011, 2016.

\bibitem{guo2012modeling}
Jianchun Guo, Yuxuan Liu, et~al.
\newblock Modeling of proppant embedment: elastic deformation and creep
  deformation.
\newblock In {\em SPE International Production and Operations Conference \&
  Exhibition}. Society of Petroleum Engineers, 2012.

\bibitem{hkj12}
H.~Hajibeygi, D.~Kavounis, and P.~Jenny.
\newblock A hierarchical fracture model for the iterative multiscale finite
  volume method.
\newblock {\em Journal of Computational Physics}, 230(24):8729--8743, 2011.

\bibitem{houwu97}
T.~Hou and X.H. Wu.
\newblock A multiscale finite element method for elliptic problems in composite
  materials and porous media.
\newblock {\em J. Comput. Phys.}, 134:169--189, 1997.

\bibitem{jenny2005adaptive}
Patrick Jenny, Seong~H Lee, and Hamdi~A Tchelepi.
\newblock Adaptive multiscale finite-volume method for multiphase flow and
  transport in porous media.
\newblock {\em Multiscale Modeling \& Simulation}, 3(1):50--64, 2005.

\bibitem{kim2011stability2}
J~Kim, HA~Tchelepi, and R~Juanes.
\newblock Stability and convergence of sequential methods for coupled flow and
  geomechanics: Drained and undrained splits.
\newblock {\em Computer Methods in Applied Mechanics and Engineering},
  200(23):2094--2116, 2011.

\bibitem{kim2011stability1}
J~Kim, HA~Tchelepi, and R~Juanes.
\newblock Stability and convergence of sequential methods for coupled flow and
  geomechanics: Fixed-stress and fixed-strain splits.
\newblock {\em Computer Methods in Applied Mechanics and Engineering},
  200(13):1591--1606, 2011.

\bibitem{kim2010sequential}
Jihoon Kim.
\newblock {\em Sequential methods for coupled geomechanics and multiphase
  flow}.
\newblock PhD thesis, Citeseer, 2010.

\bibitem{kolesov2014splitting}
AE~Kolesov, Petr~N Vabishchevich, and Maria~V Vasilyeva.
\newblock Splitting schemes for poroelasticity and thermoelasticity problems.
\newblock {\em Computers \& Mathematics with Applications}, 67(12):2185--2198,
  2014.

\bibitem{logg2009efficient}
Anders Logg.
\newblock Efficient representation of computational meshes.
\newblock {\em International Journal of Computational Science and Engineering},
  4(4):283--295, 2009.

\bibitem{logg2012automated}
Anders Logg, Kent-Andre Mardal, and Garth Wells.
\newblock {\em Automated solution of differential equations by the finite
  element method: The FEniCS book}, volume~84.
\newblock Springer Science \& Business Media, 2012.

\bibitem{lunati2006multiscale}
Ivan Lunati and Patrick Jenny.
\newblock Multiscale finite-volume method for compressible multiphase flow in
  porous media.
\newblock {\em Journal of Computational Physics}, 216(2):616--636, 2006.

\bibitem{martin2005modeling}
Vincent Martin, J{\'e}r{\^o}me Jaffr{\'e}, and Jean~E Roberts.
\newblock Modeling fractures and barriers as interfaces for flow in porous
  media.
\newblock {\em SIAM Journal on Scientific Computing}, 26(5):1667--1691, 2005.

\bibitem{olorode17}
Y.Efendiev O.M.Olorode, I.Y.Akkutlu.
\newblock A compositional model for co2 storage in deformable organic-rich
  shales.
\newblock In {\em SPE}. Society of Petroleum Engineers, 2017.

\bibitem{salimzadeh2017three}
Saeed Salimzadeh, Adriana Paluszny, and Robert~W Zimmerman.
\newblock Three-dimensional poroelastic effects during hydraulic fracturing in
  permeable rocks.
\newblock {\em International Journal of Solids and Structures}, 108:153--163,
  2017.

\bibitem{tene2016multiscale}
M~Tene, MS~Al~Kobaisi, and H~Hajibeygi.
\newblock Multiscale projection-based embedded discrete fracture modeling
  approach (f-ams-pedfm).
\newblock In {\em ECMOR XV-15th European Conference on the Mathematics of Oil
  Recovery}, 2016.

\bibitem{ctene2016algebraic}
Matei {\c{T}}ene, Mohammed~Saad Al~Kobaisi, and Hadi Hajibeygi.
\newblock Algebraic multiscale method for flow in heterogeneous porous media
  with embedded discrete fractures (f-ams).
\newblock {\em Journal of Computational Physics}, 321:819--845, 2016.

\bibitem{ctene2017projection}
Matei {\c{T}}ene, Sebastian~BM Bosma, Mohammed~Saad Al~Kobaisi, and Hadi
  Hajibeygi.
\newblock Projection-based embedded discrete fracture model (pedfm).
\newblock {\em Advances in Water Resources}, 105:205--216, 2017.

\bibitem{warren1963behavior}
JE~Warren, P~Jj Root, et~al.
\newblock The behavior of naturally fractured reservoirs.
\newblock {\em Society of Petroleum Engineers Journal}, 3(03):245--255, 1963.

\bibitem{wasaki2015permeability}
Asana Wasaki, I~Yucel Akkutlu, et~al.
\newblock Permeability of organic-rich shale.
\newblock {\em SPE Journal}, 2015.

\bibitem{weinan2007heterogeneous}
E~Weinan, Bjorn Engquist, Xiantao Li, Weiqing Ren, and Eric Vanden-Eijnden.
\newblock Heterogeneous multiscale methods: a review.
\newblock {\em Commun. Comput. Phys}, 2(3):367--450, 2007.

\bibitem{wilson1982theory}
RK~Wilson and Elias~C Aifantis.
\newblock On the theory of consolidation with double porosity.
\newblock {\em International Journal of Engineering Science}, 20(9):1009--1035,
  1982.

\bibitem{yoonspatial}
Hyun~C Yoon and Jihoon Kim.
\newblock Spatial stability for the monolithic and sequential methods with
  various space discretizations in poroelasticity.
\newblock {\em International Journal for Numerical Methods in Engineering}.

\bibitem{zhang2008sorption}
Hongbin Zhang, Jishan Liu, and D~Elsworth.
\newblock How sorption-induced matrix deformation affects gas flow in coal
  seams: a new fe model.
\newblock {\em International Journal of Rock Mechanics and Mining Sciences},
  45(8):1226--1236, 2008.

\bibitem{zhao2006fully}
Ying Zhao and Mian Chen.
\newblock Fully coupled dual-porosity model for anisotropic formations.
\newblock {\em International Journal of Rock Mechanics and Mining Sciences},
  43(7):1128--1133, 2006.

\end{thebibliography}

\end{document}